**SINIR KOŞULUNDA SPEKTRAL PARAMETRE
BULUNAN GEÇ KALAN ARGÜMENTLİ
SÜREKLİ OLMAYAN SINIR-DEĞER
PROBLEMİNİN ÖZDEĞER VE
ÖZFONKSİYONLARININ
ASİMPTOTİK İFADELERİ**

**Erdoğan ŞEN
Yüksek Lisans Tezi
Matematik Anabilim Dalı
Danışman: Doç.Dr. Azad Bayramov
2010**

T.C.

NAMIK KEMAL ÜNİVERSİTESİ

FEN BİLİMLERİ ENSTİTÜSÜ

YÜKSEK LİSANS TEZİ

**SINIR KOŞULUNDA SPEKTRAL PARAMETRE BULUNAN GEÇ KALAN ARGÜMENTLİ SÜREKLİ OLMAYAN SINIR-DEĞER PROBLEMİNİN ÖZDEĞER VE ÖZFONKSİYONLARININ ASİMPTOTİK İFADELERİ**

Erdoğan ŞEN

**MATEMATİK ANABİLİM DALI**

**DANIŞMAN: Doç. Dr. Azad BAYRAMOV**

**TEKİRDAĞ-2010**

**Her hakkı saklıdır**

Doç. Dr. Azad BAYRAMOV'un .danışmanlığında, Erdoğan ŞEN tarafından hazırlanan bu çalışma aşağıdaki jüri tarafından Matematik Anabilim Dalı'nda yüksek lisans tezi olarak kabul edilmiştir.

Juri Başkanı : Prof. Dr. Rıfat MİRKASIM　　　　　　　　*İmza* :

Üye : Doç. Dr. Azad BAYRAMOV (Danışman)　　　　　*İmza* :

Üye : Yrd. Doç. Dr. Dilek Çiftçi KAZICI　　　　　　　*İmza* :

Fen Bilimleri Enstitüsü Yönetim Kurulunun ………………. tarih ve ………………. sayılı kararıyla onaylanmıştır.

Doç. Dr. Fatih KONUKCU
**Enstitü Müdürü**

# ÖZET

Yüksek Lisans Tezi

SINIR KOŞULUNDA SPEKTRAL PARAMETRE BULUNAN GEÇ KALAN ARGÜMENTLİ SÜREKLİ OLMAYAN SINIR-DEĞER PROBLEMİNİN ÖZDEĞER VE ÖZFONKSİYONLARININ ASİMPTOTİK İFADELERİ

Erdoğan ŞEN

Namık Kemal Üniversitesi
Fen Bilimleri Enstitüsü
Matematik Anabilim Dalı

Danışman : Doç. Dr. Azad BAYRAMOV

Bu Araştırmanın amacı sınır koşulunda spektral parametre bulunan geç kalan argümentli sürekli olmayan sınır -değer probleminin özdeğer ve özfonksiyonlarının asimptotik ifadelerinin bulunmasıdır.

Otomatik Kontrol kuramında, öz-titreşim sistemleri kuramında, roket motorlarının ateşlenmesi ile İlgili çalışmalarda, ekonominin, biyofiziğin ve daha başka alanların birçok problemlerinde geç kalan argümentli diferansiyel denklemlerin uygulamalarına rastlanılır. Bu alanlardaki problemler geç kalan argümentli diferansiyel denklemlere indirgenerek çözülür. Bu çalışmada da sınır koşulunda spektral parametre bulunan geç kalan argümentli sürekli olmayan sınır-değer problemi incelenmiş, yukarıda bahsedilen alanlarda kullanılmak üzere özdeğer ve özfonksiyonlar için asimptotik formüller elde edilmiştir.

**Anahtar kelimeler:** Geç kalan argümentli diferansiyel denklem, geçiş koşulları, özdeğer ve özfonksiyonların asimptotikleri, spektral parametre

**2010, 27 sayfa**

# ABSTRACT

MSc. Thesis

ASYMPTOTIC EXPRESSIONS OF EIGENVALUES AND EIGENFUNCTIONS OF A DISCONTINUOUS BOUNDARY-VALUE PROBLEM WITH RETARDED ARGUMENT WHICH CONTAINS A SPECTRAL PARAMETER IN THE BOUNDARY CONDITION.

Erdoğan ŞEN

Namık Kemal University
Graduate School of Natural and Applied Sciences
Department of Mathematics

Supervisor : Assoc. Prof. Dr. Azad BAYRAMOV

The aim of this study is to find asymptotic expressions of eigenvalues and eigenfunctions of a discontinuous boundary-value problem with retarded argument which contains a spectral parameter in the boundary condition.

Applications of differential equations with retarded argument can be encountered in the theory of self-oscillatory systems, in the study of problems connected with combustion in rocket engines, in a number of problems in economics, biophysics, and many other fields. The problems in these areas can be solved reducing differential equations with retarded argument. In this study discontinuous boundary-value problem with retarded argument which contains a spectral parameter in the boundary condition were investigated and asymptotic formulas were obtained for eigenvalues and eigenfunctions for using areas which mentioned above.

**2010, 27 pages**

# TEŞEKKÜR

Bu çalışmada bana destek olan ve emeği geçen danışman hocam Sayın Doç. Dr. Azad BAYRAMOV'a teşekkür ederim.

# SİMGELER VE KISALTMALAR DİZİNİ

| | |
|---|---|
| $O$ | : Büyük $O$ notasyonu |
| $y'$ | : $y$ fonksiyonunun birinci mertebeden türevi |
| $y''$ | : $y$ fonksiyonunun ikinci mertebeden türevi |
| $\Delta$ | : Geç kalan argüment |

## İÇİNDEKİLER

## 1. GİRİŞ

Geç kalan argümentli diferansiyel denklemler ilk olarak XVIII. yüzyılda Euler probleminin çözümüyle bağlantılı olarak ortaya çıkmıştır. Geç kalan argümentli diferansiyel denklemler adi diferansiyel denklemler teorisinin güncel konularından biridir. Özellikle geç kalan argümentli diferansiyel denklemlerin spektral analizine gittikçe artan bir ilgi vardır. Matematiksel fiziğin bir çok problemi geç kalan argümentli diferansiyel denklemlere indirgenerek çözülür. Sınır koşulunda spektral parametre bulunan ikinci mertebeden adi diferansiyel operatörler için sınır-değer problemi Fulton (1972), Kerimov ve Mamedov (1999), Muktarov ve ark. (2003) ve Tikhonov (1972) tarafından çalışılmıştır. İkinci mertebeden geç kalan argümentli diferansiyel denklemler için Sturm-Liouville tipli sınır-değer probleminin özdeğer ve özfonksiyonlarının asimptotik özellikleri ise Norkin (1958 ve 1972), Bellman ve Cook (1963), Demidenko ve Likhoshvai (2005), Bayramov ve ark. (2007), Bayramoğlu ve ark. (2002) tarafından çalışılmıştır. Ayrıca sınır koşulunda spektral parametre bulunan Sturm-Liouville tipli geç kalan argümentli diferansiyel denklem için özdeğer ve özfonksiyonların asimptotik formülleri Bayramoğlu ve ark. (2002) tarafından elde edilmiştir. Bu çalışmada da $[0,\pi]$ aralığında $x = \pi/2$ süreksizlik noktasını içeren,

$$p(x)y''(x) + q(x)y(x - \Delta(x)) + \lambda y(x) = 0 \quad (1.1)$$

diferansiyel denklemini

$$a_1 y(0) + a_2 y'(0) = 0, \quad (1.2)$$

$$y'(\pi) + d\lambda y(\pi) = 0 \quad (1.3)$$

sınır koşulları ve

$$\gamma_1 y(\frac{\pi}{2} - 0) = \delta_1 y(\frac{\pi}{2} + 0), \quad (1.4)$$

$$\gamma_2 y'(\frac{\pi}{2} - 0) = \delta_2 y'(\frac{\pi}{2} + 0) \quad (1.5)$$

taşıma koşulları ile göz önüne alacağız.

Bu denklemde $p(x)$, $q(x)$ ve $\Delta(x) \geq 0$, $[0,\frac{\pi}{2}) \bigcup (\frac{\pi}{2},\pi]$ aralığında sürekli ve reel değerli fonksiyonlardır. Burada $p(x)$

$$p(x) = \begin{cases} p_1^2, & x \in [0, \frac{\pi}{2}), \\ p_2^2, & x \in (\frac{\pi}{2}, \pi] \end{cases}$$

şeklinde tanımlıdır ve şu koşullar sağlanır: $q(\frac{\pi}{2} \pm 0) = \lim_{x \to \frac{\pi}{2} \pm 0} q(x)$, $\Delta(\frac{\pi}{2} \pm 0) = \lim_{x \to \frac{\pi}{2} \pm 0} \Delta(x)$ sonlu limitleri mevcuttur, $x \in [0, \frac{\pi}{2})$ için $x - \Delta(x) \geq 0$, $x \in (\frac{\pi}{2}, \pi]$ için $x - \Delta(x) \geq \frac{\pi}{2}$. $\lambda$ reel spektral parametre, $p_1, p_2, \gamma_1, \gamma_2, \delta_1, \delta_2, a_1, a_2, d$ reel sayılar, $|a_1| + |a_2| \neq 0$ ve $i = 1, 2$ için $|\gamma_i| + |\delta_i| \neq 0$. Ayrıca varsayalım ki $\gamma_1 \delta_2 p_1 = \gamma_2 \delta_1 p_2$ eşitliği sağlanır.

## 2. KURAMSAL TEMELLER

İkinci mertebeden sapan argümentli diferansiyel denklemler

$$F\left(t, x(t),...,x^{(m_0)}(t), x(t-\Delta_1(t)),...,x^{(m_1)}(t-\Delta_1(t)),...,x(t-\Delta_n(t)),...,x^{(m_n)}(t-\Delta_n(t))\right) = 0 \quad (2.1)$$

şeklindedir. Burada $i = 1,...,n$ için $\Delta_i(t) \geq 0$ sürekli ve $\max_{0 \leq i \leq n} m_i = 2$ olarak verilir. $x^{(m_i)}(t-\Delta_i(t))$ ile $x(z)$ fonksiyonunun $z = t - \Delta_i(t)$ noktasındaki türevi kastedilmektedir. $A$ verilen başlangıç noktası olsun. Her $\Delta_i(t)$ sapması $A$ noktasını içeren bir $E_A^{(i)}$ başlangıç kümesi tanımlar ve $t \geq A$ için $t - \Delta_i(t) < A$ dır. $E_A = \bigcup_{i=1}^{n} E_A^{(i)}$ ve $\mu = \max_{1 \leq i \leq n} m_i$ olsun. $E_A$ üzerinde $\mu$-kere türevlenebilen bir $\Phi(t)$ başlangıç fonksiyonu tayin edelim.

$x_A^{(j)} = \Phi^{(j)}(A), j = 0,...,\mu$ olsun. $\mu = 0$ ise ek olarak $x_A^{(1)}$ sayısını tayin edelim. Eğer $A$, $E_A$ kümesinin izole edilmiş bir noktası ise $x_A^{(0)}$ ve $x_A^{(1)}$ keyfi olarak seçilir.

(2.1) denklemi için başlangıç değer problemi $[A, B), B \leq +\infty$ aralığında

$$x(A) = x_A^{(0)}, \; x'(A) = x_A^{(1)},$$
$$x^{(j)}\left(t - \Delta_i(t)\right) \equiv \Phi^{(j)}\left(t - \Delta_i(t)\right), \; t - \Delta_i(t) < A \; ise \quad (2.2)$$

koşullarını sağlayan $x(t)$ çözümünü bulma problemidir.

Sapan argümentli diferansiyel denklemlerin doğal bir sınıflandırılması G. A. Kamenskii tarafından yapılmıştır. (2.1) denklemi $x^{(m_0)}(t)$ için çözülürse

$$x^{(m_0)}(t) = f\left(t, x(t), x^{(m_0-1)}(t), x(t-\Delta_1(t)),...,x^{(m_1)}(t-\Delta_1(t))... \right.$$
$$\left. ...,x(t-\Delta_n(t)),...,x^{(m_n)}(t-\Delta_n(t))\right) \quad (2.3)$$

olur. $\lambda = m_0 - \mu$ olsun. $\lambda > 0$ için denklemler geç kalan argümentli denklemler; $\lambda = 0$ için nötral tipli denklemler ve $\lambda < 0$ için ileri tipli denklemler olarak adlandırılır.

(2.3) denkleminde $\lambda > 0$ ve $f$ fonksiyonu $t$ nin dışındaki tüm argümentlere göre lineerse ikinci mertebeden geç kalan argümentli diferansiyel denklem elde ederiz:

$$x''(t) = \sum_{i=0}^{n} \left( a_i(t) x(t - \Delta_i(t)) + b_i(t) x'(t - \Delta_i(t)) \right) + c(t). \tag{2.4}$$

(2.4) denkleminde $x(t) = y_1(t)$ ve $x'(t) = y_2(t)$ olsun. (2.4) denklemini birinci mertebeden

$$y_1'(t) = y_2(t)$$
$$y_2'(t) = \sum_{i=0}^{n} \left( a_i(t) y_1(t - \Delta_i(t)) + b_i(t) y_2(t - \Delta_i(t)) \right) + c(t)$$

sistemi ile yer değiştirelim. Uygunluk için daha genel

$$y_k'(t) = \sum_{j=1}^{2} \sum_{i=0}^{n} a_{ji}^{k}(t) y_j(t - \Delta_i(t)) + c_k(t), \quad k = 1, 2. \tag{2.5}$$

sistemini göz önüne alalım. (2.5) sistemi dağıtılmış gecikmelerle birlikte

$$y_k'(t) = \sum_{j=1}^{m} \int_{0}^{\infty} y_j(t - s) \, dr_j^k(t, s) + c_k(t), \quad k = 1, 2, ..., m \tag{2.6}$$

sisteminin özel bir biçimidir. Burada integral Stieltjes anlamında alınmıştır.

Myskis (1951) $r_j^k(t, s)$ üzerinde bazı kısıtlamalara giderek (2.6) sisteminin başlangıç –değer probleminin çözümleri için varlık ve teklik teoremleri oluşturmuştur. Dahası bu çözümler başlangıç verilerine sürekli olarak bağlıdır. Amacımıza uygun olarak bu teoremleri (2.5) sisteminin

$$y_k'(t) = \sum_{j=1}^{2} \left( a_{j0}^{k}(t) y_j(t) + a_{j1}^{k}(t) y_j(t - \Delta(t)) \right) + c_k(t), \quad k = 1, 2. \tag{2.7}$$

formuna uygun olarak yazmak yeterlidir.

(2.7) eşitliğinin $[A, B)$, $B \leq +\infty$ aralığında tanımlı

$$y_k(A) = y_A^{(k)} = \Phi_k(A),$$
$$y_k(t - \Delta(t)) \equiv \Phi_k(t - \Delta(t)), \quad t - \Delta(t) < A \tag{2.8}$$

koşullarını sağlayan $y_1(t)$, $y_2(t)$ çözümlerini araştıracağız. Burada $k = 1, 2$ ve $\Phi_k(t)$, $E_A$ üzerinde tanımlı başlangıç fonksiyonu olsun. $a_{j1}^{k}(t) \equiv 0$ ($j = 1$ veya 2) ise, $y_A^{(j)}$ ($j = 1$ veya

2) sayısı keyfi olarak atanır; eğer $A$ noktası $E_A$ da izole edilmiş bir nokta ise $y_A^{(1)}$ ve $y_A^{(2)}$ sayılarının her ikisi de keyfi olarak atanır.

(2.7) sistemi ile birlikte

$$y_k(t) = y_A^{(k)} + \int_A^t \left\{ \sum_{j=1}^2 \left( a_{j0}^{\ k}(\tau) y_j(\tau) + a_{j1}^{\ k}(\tau) \right) y_j(\tau - \Delta(t)) \right\} d\tau \qquad (2.9)$$
$$+ \int_A^t c_k(\tau) d\tau, \qquad k=1,2; \qquad A \leq t < B$$

integral denklemler sistemini göz önüne alalım. (2.8) e eşdeğer olarak da

$$y_j(\tau - \Delta(\tau)) \equiv \Phi_j(\tau - \Delta(\tau)), \quad \tau - \Delta(\tau) < A, \qquad j=1,2. \qquad (2.10)$$

koşulunu göz önüne alalım.

Bundan böyle $a_{j0}^{\ k}(t), a_{j1}^{\ k}(t), c_k(t)$ ($j,k=1,2$) fonksiyonlarının ve $\Delta(t) \geq 0$ ın $[A,B)$ üzerinde ve $\Phi_k(t)(k=1,2)$ başlangıç fonksiyonlarının $E_A$ üzerinde sürekli olduklarını varsayacağız.

**Lemma 2.1** (2.7) sisteminin (2.8) koşullarını sağlayan bir çözümü, (2.9) integral denklemler sisteminin (2.10) koşulunu sağlayan sürekli bir çözümüdür. Tersine (2.9) sisteminin (2.10) koşulunu sağlayan sürekli bir çözümü (2.7) sisteminin (2.8) koşullarını sağlayan bir çözümüdür.

**Teorem 2.1** $\Phi = \max_{1 \leq k \leq 2} \sup_{E_A} |\Phi_k(t)| < \infty$ olsun. O halde (2.7) sisteminin $[A,B)$ aralığında (2.8) başlangıç koşullarını sağlayan tek bir çözümü vardır.

## 3. MATERYAL ve YÖNTEM

Bu çalışmada öncelikle (1.1)-(1.5) probleminin çözümü integral denklemler cinsinden yazılmıştır. Daha sonra problemin özdeğerlerinin sayısı ve yapısı belirlenerek özdeğer ve özfonksiyonlar için asimptotik formüller elde edilmiştir.

Teoremlerin ispatlanmasında ise Rolle teoremi, kısmi integrasyon ve iterasyon tekniği kullanılmıştır.

**Rolle Teoremi.** $f(x)$ fonksiyonu $[a,b]$ aralığında sürekli, $f'(x)$ türevi $(a,b)$ açık aralığında mevcut ve $f(a) = f(b)$ olsun. O halde en az bir $c \in (a,b)$ için $f'(c) = 0$ olur.

**Teorem (Kısmi integrasyon).** $u(x)$ ve $v(x)$ fonksiyonları $[a,b]$ kapalı aralığında diferansiyellenebilen fonksiyonlar ise $\int_a^b u(x)v'(x)dx = u(b)v(b) - u(a)v(a) - \int_a^b u'(x)v(x)dx$ olur.

**İterasyon.** Çekirdek olarak adlandırılan $K(t,s)$ ve $f(t)$ fonksiyonları bilinen, $y(t)$ bilinmeyen fonksiyon, $\lambda$ ise herhangi bir sayısal parametre olmak üzere

$$y(t) = \lambda \int_a^b K(t,s)y(s)ds + f(t), \quad t \in [a,b] \tag{3.1}$$

lineer Fredholm integral denklemini göz önüne alalım. Burada $K(t,s)$ fonksiyonu $G = \{(t,s) \in R^2 : a \leq t \leq b, a \leq s \leq b\}$ üzerinde, $f(t)$ ise $[a,b]$ üzerinde süreklidir. O halde (3.1) denkleminin $[a,b]$ aralığı üzerinde sürekli bir $y^*(t)$ çözümü vardır ve her başlangıç $y_0(t) \in C[a,b]$ fonksiyonu için terimleri

$$y_n(t) = \lambda \int_a^b K(t,s)y_{n-1}(s)ds + f(t), \quad n = 1, 2, \ldots$$

şeklinde tanımlanan $y_n$ dizisi $y^*$ fonksiyonuna $[a,b]$ aralığında düzgün yakınsaktır.

## 4. ARAŞTIRMA BULGULARI

$w_1(x,\lambda)$, (1.1) denkleminin $\left[0, \dfrac{\pi}{2}\right]$ aralığında

$$w_1(0,\lambda) = a_2, \quad w_1'(0,\lambda) = -a_1 \tag{4.1}$$

başlangıç koşullarını sağlayan bir çözümü olsun. (4.1) başlangıç koşulları (1.1) denkleminin tek bir çözümünü tanımlar (Norkin 1972). Yukarıdaki durumu belirterek $\left[\dfrac{\pi}{2}, \pi\right]$ aralığında (1.1) denkleminin $w_2(x,\lambda)$ çözümünü aşağıdaki başlangıç koşullarında $w_1(x,\lambda)$ çözümü yardımıyla şöyle tanımlayacağız:

$$\begin{aligned} w_2(\dfrac{\pi}{2},\lambda) &= \gamma_1 \delta_1^{-1} w_1(\dfrac{\pi}{2},\lambda), \\ w_2'(\dfrac{\pi}{2},\lambda) &= \gamma_2 \delta_2^{-1} w_1'(\dfrac{\pi}{2},\lambda). \end{aligned} \tag{4.2}$$

(4.2) başlangıç koşulları (1.1) denkleminin $\left[\dfrac{\pi}{2}, \pi\right]$ aralığında tek bir çözümünü tanımlar. Sonuç olarak $w(x,\lambda)$ fonksiyonu $\left[0, \dfrac{\pi}{2}\right) \bigcup \left(\dfrac{\pi}{2}, \pi\right]$ aralığında (1.1) denkleminin (1.2),(1.4) ve (1.5) sınır koşullarını sağlayan çözümü olarak aşağıdaki eşitlik ile tanımlıdır.

$$w(x,\lambda) = \begin{cases} w_1(x,\lambda), & x \in \left[0, \dfrac{\pi}{2}\right), \\ w_2(x,\lambda), & x \in \left(\dfrac{\pi}{2}, \pi\right]. \end{cases}$$

**Lemma 4.1** $w(x,\lambda)$, (1.1) denkleminin bir çözümü ve $\lambda > 0$ olsun. O halde aşağıdaki integral denklemler sağlanır.

$$w_1(x,\lambda) = a_2 \cos\dfrac{s}{p_1}x - \dfrac{a_1 p_1}{s}\sin\dfrac{s}{p_1}x - \dfrac{1}{s}\int_0^x \dfrac{q(\tau)}{p_1}\sin\dfrac{s}{p_1}(x-\tau) w_1(\tau - \Delta(\tau),\lambda)\,d\tau, \tag{4.3}$$

$$w_2(x,\lambda) = \frac{\gamma_1}{\delta_1} w_1(\frac{\pi}{2},\lambda) \cos\frac{s}{p_2}(x-\frac{\pi}{2}) + \frac{\gamma_2 p_2 w_1'(\frac{\pi}{2},\lambda)}{s\delta_2} \sin\frac{s}{p_2}(x-\frac{\pi}{2})$$
$$-\frac{1}{s}\int_{\pi/2}^{x} \frac{q(\tau)}{p_2} \sin\frac{s}{p_2}(x-\tau) w_2(\tau-\Delta(\tau),\lambda) d\tau \quad (s=\sqrt{\lambda}, \lambda>0). \tag{4.4}$$

**İspat.** Bu lemmayı ispatlamak için (4.3) ve (4.4) de sırasıyla $-\frac{q(\tau)}{p_1^2} w_1(\tau-\Delta(\tau),\lambda)$ yerine

$-\frac{s^2}{p_1^2} w_1(\tau,\lambda) - w_1''(\tau,\lambda)$ ve $-\frac{q(\tau)}{p_2^2} w_2(\tau-\Delta(\tau),\lambda)$ yerine $-\frac{s^2}{p_2^2} w_2(\tau,\lambda) - w_2''(\tau,\lambda)$

konularak ve iki kere kısmi integrasyon uygulanarak bulunur. ∎

**Teorem 4.1** (1.1)-(1.5) sınır-değer problemi sadece basit özdeğerlere sahip olabilir.

**İspat.** $\tilde{\lambda}$, (1.1)-(1.5) sınır-değer probleminin bir özdeğeri ve

$$\tilde{u}(x,\tilde{\lambda}) = \begin{cases} \tilde{u}_1(x,\tilde{\lambda}), & x \in \left[0,\frac{\pi}{2}\right), \\ \tilde{u}_2(x,\tilde{\lambda}), & x \in \left(\frac{\pi}{2},\pi\right] \end{cases}$$

bu özdeğere karşı gelen özfonksiyon olsun. O halde (1.2) ve (4.1) den

$$W\left[\tilde{u}_1(0,\tilde{\lambda}), w_1(0,\tilde{\lambda})\right] = \begin{vmatrix} \tilde{u}_1(0,\tilde{\lambda}) & a_2 \\ \tilde{u}_1'(0,\tilde{\lambda}) & -a_1 \end{vmatrix} = 0$$

elde edilir. Böylece Wronskian sıfıra eşit olduğundan $\tilde{u}_1(x,\tilde{\lambda})$ ve $w_1(x,\tilde{\lambda})$, $\left[0,\frac{\pi}{2}\right]$ aralığında lineer bağımlıdır (Norkin 1972). Benzer şekilde $\tilde{u}_2(x,\tilde{\lambda})$ ve $w_2(x,\tilde{\lambda})$ nün de $\left[\frac{\pi}{2},\pi\right]$ aralığında lineer bağımlı olduğu gösterilebilir. Bu yüzden bazı $K_1 \neq 0$ ve $K_2 \neq 0$ için

$$\tilde{u}_i(x,\tilde{\lambda}) = K_i w_i(x,\tilde{\lambda}) \quad (i=1,2) \tag{4.5}$$

elde edilir. $K_1 = K_2$ olduğunu göstermeliyiz. Varsayalım ki $K_1 \neq K_2$. (1.4) ve (4.5) eşitliklerinden

$$\gamma_1 \tilde{u}(\frac{\pi}{2}-0,\tilde{\lambda}) - \delta_1 \tilde{u}(\frac{\pi}{2}+0,\tilde{\lambda}) = \gamma_1 \tilde{u}_1(\frac{\pi}{2},\tilde{\lambda}) - \delta_1 \tilde{u}_2(\frac{\pi}{2},\tilde{\lambda}) = \gamma_1 K_1 w_1(\frac{\pi}{2},\tilde{\lambda}) - \delta_1 K_2 w_2(\frac{\pi}{2},\tilde{\lambda})$$

$$= \gamma_1 K_1 \delta_1 \gamma_1^{-1} w_2(\frac{\pi}{2},\tilde{\lambda}) - \delta_1 K_2 w_2(\frac{\pi}{2},\tilde{\lambda}) = K_1 \delta_1 w_2(\frac{\pi}{2},\tilde{\lambda}) - K_2 \delta_1 w_2(\frac{\pi}{2},\tilde{\lambda})$$

$$= \delta_1 (K_1 - K_2) w_2(\frac{\pi}{2},\tilde{\lambda}) = 0.$$

$\delta_1(K_1 - K_2) \neq 0$ olduğundan

$$w_2(\frac{\pi}{2},\tilde{\lambda}) = 0 \tag{4.6}$$

elde edilir. Benzer şekilde (1.5) den

$$w_2'(\frac{\pi}{2},\tilde{\lambda}) = 0 \tag{4.7}$$

elde edilir. $w_2(x,\tilde{\lambda})$, $\left[\frac{\pi}{2},\pi\right]$ aralığında (1.1) denkleminin (4.6) ve (4.7) başlangıç koşullarını sağlayan bir çözümü olduğundan $\left[\frac{\pi}{2},\pi\right]$ aralığında $w_2(x,\tilde{\lambda}) = 0$ olarak bulunur. Ayrıca,

$$w_1(\frac{\pi}{2},\tilde{\lambda}) = w_1'(\frac{\pi}{2},\tilde{\lambda}) = 0.$$

eşitliği (4.2), (4.6) ve (4.7) kullanarak elde edilir.

$w_2(x,\lambda)$ hakkında elde edilen sonuçlardan $\left[0,\frac{\pi}{2}\right]$ aralığında $w_1(x,\lambda) = 0$ olarak bulunur. Dolayısıyla $\left[0,\frac{\pi}{2}\right) \cup \left(\frac{\pi}{2},\pi\right]$ aralığında $w(x,\tilde{\lambda})$ özdeş olarak sıfıra eşittir. Bu da (4.1) ile çelişir ki bu da ispatı tamamlar. ∎

$w(x,\lambda)$, (1.1) denkleminin (1.2), (1.4) ve (1.5) koşullarını sağlayan aşikar olmayan bir çözümüdür. $w(x,\lambda)$ yı (1.3) de yerine koyarsak,

$$F(\lambda) = w'(\pi,\lambda) + d\lambda w(\pi,\lambda) = 0. \tag{4.8}$$

karakteristik denklemini elde ederiz.

Teorem (4.1) den (1.1)-(1.5) sınır-değer probleminin özdeğerleri ile (4.8) denkleminin reel kökleri aynıdır.

$$q_1 = \frac{1}{p_1} \int_0^{\pi/2} |q(\tau)| d\tau \text{ ve } q_2 = \frac{1}{p_2} \int_{\pi/2}^{\pi} q(\tau) d\tau \text{ olsun.}$$

**Lemma 4.2** (1) $\lambda \geq 4q_1^2$ olsun. O halde (4.3) denkleminin $w_1(x,\lambda)$ çözümü için aşağıdaki eşitsizlik doğrudur:

$$|w_1(x,\lambda)| \leq \frac{1}{q_1}\sqrt{4q_1^2 a_2^2 + p_1^2 a_1^2}, \qquad x \in \left[0, \frac{\pi}{2}\right]. \tag{4.9}$$

(2) $\lambda \geq \max\{4q_1^2, 4q_2^2\}$ olsun. O halde (4.4) denkleminin $w_2(x,\lambda)$ çözümü için aşağıdaki eşitsizlik doğrudur:

$$|w_2(x,\lambda)| \leq \frac{2}{q_1}\sqrt{4q_1^2 a_2^2 + p_1^2 a_1^2}\left\{\frac{|\gamma_1|}{|\delta_1|} + \left|\frac{p_2 \gamma_2}{p_1 \delta_1}\right|\right\}, \qquad x \in \left[\frac{\pi}{2}, \pi\right]. \tag{4.10}$$

**İspat.** $B_{1\lambda} = \max_{\left[0,\frac{\pi}{2}\right]}|w_1(x,\lambda)|$ olsun. O halde (4.3) den her $\lambda > 0$ için aşağıdaki eşitsizlik doğrudur:

$$B_{1\lambda} \leq \sqrt{a_2^2 + \frac{p_1^2 a_1^2}{s^2}} + \frac{1}{s}B_{1\lambda}q_1.$$

Eğer $s \geq 2q_1$ ise o halde

$$B_{1\lambda} \leq \sqrt{\frac{s^2 a_2^2 + p_1^2 a_1^2}{s^2}} + \frac{1}{s}B_{1\lambda}q_1 = \frac{1}{s}\sqrt{s^2 a_2^2 + p_1^2 a_1^2} + \frac{B_{1\lambda}}{2}.$$

Buradan da

$$B_{1\lambda} - \frac{B_{1\lambda}}{2} = \frac{B_{1\lambda}}{2} \leq \frac{1}{s}\sqrt{s^2 a_2^2 + p_1^2 a_1^2} \leq \frac{1}{2q_1}\sqrt{4q_1^2 a_2^2 + p_1^2 a_1^2}$$

$$\Rightarrow B_{1\lambda} \leq \frac{1}{q_1}\sqrt{4q_1^2 a_2^2 + p_1^2 a_1^2}.$$

Bulunur ki böylece (4.9) elde edilir.

(4.3) ün $x$ e göre türevini alırsak şunu elde ederiz:

$$w_1'(x,\lambda) = -\frac{s}{p_1}a_2 \sin\frac{s}{p_1}x - a_1 \cos\frac{s}{p_1}x - \frac{1}{p_1^2}\int_0^x q(\tau)\cos\frac{s}{p_1}(x-\tau)w_1(\tau-\Delta(\tau))d\tau. \tag{4.11}$$

$\left|-\frac{s}{p_1}a_2 \sin\frac{s}{p_1}x - a_1 \cos\frac{s}{p_1}x\right| \leq \sqrt{\frac{s^2}{p_1^2}a_2^2 + a_1^2}$ olduğunu görmek kolaydır.

$$\left|\frac{1}{p_1^2}\int_0^x q(\tau)\cos\frac{s}{p_1}(x-\tau)w_1(\tau-\Delta(\tau))d\tau\right|$$

$$\leq \frac{1}{p_1}\left(\frac{1}{p_1}\int_0^x q(\tau)d\tau\right)\left(\frac{1}{q_1}\sqrt{4q_1^2 a_2^2 + p_1^2 a_1^2}\right)\max_{\tau\in[0,x]}\left\{\cos\frac{s}{p_1}(x-\tau)\right\}$$

$$\leq \frac{1}{p_1}\cdot q_1 \cdot\left(\frac{1}{q_1}\sqrt{4q_1^2 a_2^2 + p_1^2 a_1^2}\right)\cdot 1 = \frac{1}{p_1}\sqrt{4q_1^2 a_2^2 + p_1^2 a_1^2}.$$

(4.9) ve (4.11) den $s \geq 2q_1$ için aşağıdaki eşitsizlik doğrudur:

$$|w_1'(x,\lambda)| \leq \sqrt{\frac{s^2}{p_1^2}a_2^2 + a_1^2} + \frac{1}{p_1}\sqrt{4q_1^2 a_2^2 + p_1^2 a_1^2}.$$

Bundan dolayı $s \geq 2q_1$ ve $x \in \left[0, \frac{\pi}{2}\right]$ için

$$|w_1'(x,\lambda)| \leq \frac{1}{p_1}\sqrt{4q_1^2 a_2^2 + p_1^2 a_1^2} + \frac{1}{p_1}\sqrt{4q_1^2 a_2^2 + p_1^2 a_1^2} = \frac{2}{p_1}\sqrt{4q_1^2 a_2^2 + p_1^2 a_1^2}$$

$$\leq \frac{s}{p_1 q_1}\sqrt{4q_1^2 a_2^2 + p_1^2 a_1^2}.$$

Buradan da

$$\frac{|w_1'(x,\lambda)|}{s} \leq \frac{1}{p_1 q_1}\sqrt{4q_1^2 a_2^2 + p_1^2 a_1^2} \tag{4.12}$$

elde edilir.

$B_{2\lambda} = \max_{\left[\frac{\pi}{2},\pi\right]}|w_2(x,\lambda)|$ olsun. O halde (4.4), (4.9) ve (4.12) den $s \geq 2q_1$ için

$$|w_2(x,\lambda)| \leq B_{2\lambda} \leq \left|\frac{\gamma_1}{\delta_1}\right|\left|w_1(\frac{\pi}{2},\lambda)\right| + |p_2|\left|\frac{\gamma_2}{\delta_2}\right|\frac{\left|w_1'(\frac{\pi}{2},\lambda)\right|}{s} + \frac{1}{s}B_{2\lambda}\int_{\pi/2}^{x}\frac{q(\tau)}{p_2}d\tau$$

$$\leq \frac{1}{q_1}\left|\frac{\gamma_1}{\delta_1}\right|\sqrt{4q_1^2 a_2^2 + p_1^2 a_1^2} + |p_2|\left|\frac{\gamma_2}{\delta_2}\right|\frac{1}{|p_1 q_1|}\sqrt{4q_1^2 a_2^2 + p_1^2 a_1^2} + \frac{1}{2q_2}B_{2\lambda}q_2$$

$$\Rightarrow B_{2\lambda} - \frac{B_{2\lambda}}{2} = \frac{B_{2\lambda}}{2} \leq \frac{1}{q_1}\sqrt{4q_1^2 a_2^2 + a_1^2 p_1^2}\left\{\left|\frac{\gamma_1}{\delta_1}\right| + \left|\frac{p_2 \gamma_2}{p_1 \delta_1}\right|\right\}.$$

Bundan dolayı $\lambda \geq \max\{4q_1^2, 4q_2^2\}$ için (4.10) sağlanır. ∎

**Teorem 4.2** (1.1)-(1.5) problemi sonsuz sayıda pozitif özdeğerlere sahiptir.

**İspat.** (4.4) ün $x$ e göre türevini alırsak

$$w_2'(x,\lambda) = -\frac{s}{p_2}\frac{\gamma_1}{\delta_1}w_1'\left(\frac{\pi}{2},\lambda\right)\sin\frac{s}{p_2}(x-\frac{\pi}{2}) + \frac{\gamma_2}{\delta_2}w_1'\left(\frac{\pi}{2},\lambda\right)\cos\frac{s}{p_2}(x-\frac{\pi}{2})$$

$$-\frac{1}{p_2^2}\int_{\pi/2}^{x}q(\tau)\cos\frac{s}{p_2}(x-\frac{\pi}{2})w_2(\tau-\Delta(\tau),\lambda)d\tau. \tag{4.13}$$

(4.3), (4.4), (4.6), (4.9) ve (4.11) den

$$F(\lambda) = -\frac{s}{p_2}\frac{\gamma_1}{\delta_1}\left(a_2 \cos\frac{s\pi}{p_1 2} - \frac{a_1}{s}p_1 \sin\frac{s\pi}{p_1 2} - \frac{1}{sp_1}\int_0^{\pi/2} q(\tau)\sin\frac{s}{p_1}\left(\frac{\pi}{2}-\tau\right)w_1(\tau-\Delta(\tau),\lambda)d\tau\right)$$

$$\times \sin\frac{s\pi}{p_2 2}$$

$$+\frac{\gamma_2}{\delta_2}\left(-\frac{s}{p_1}a_2 \sin\frac{s\pi}{p_1 2} - a_1 \cos\frac{s\pi}{p_1 2} - \frac{1}{p_1^2}\int_0^{\pi/2} q(\tau)\cos\frac{s}{p_1}\left(\frac{\pi}{2}-\tau\right)w_1(\tau-\Delta(\tau),\lambda)\right)$$

$$\times \cos\frac{s\pi}{p_2 2}$$

$$-\frac{1}{p_2^2}\int_{\pi/2}^{\pi} q(\tau)\cos\frac{s}{p_2}(\pi-\tau)w_2(\tau-\Delta(\tau),\lambda)d\tau$$

$$+\lambda d\left(\frac{\gamma_1}{\delta_1}\left[a_2 \cos\frac{s\pi}{p_1 2} - \frac{a_1}{s}p_1 \sin\frac{s\pi}{p_1 2} - \frac{1}{sp_1}\int_0^{\pi/2} q(\tau)\sin\frac{s}{p_1}\left(\frac{\pi}{2}-\tau\right)w_1(\tau-\Delta(\tau),\lambda)d\tau\right]\right.$$

$$\times \cos\frac{s\pi}{p_2 2}$$

$$+\frac{\gamma_2 p_2}{\delta_2 s}\left[-\frac{s}{p_1}a_2 \sin\frac{s\pi}{p_1 2} - a_1 \cos\frac{s\pi}{p_1 2} - \frac{1}{p_1^2}\int_0^{\pi/2} q(\tau)\cos\frac{s}{p_1}\left(\frac{\pi}{2}-\tau\right)w_1(\tau-\Delta(\tau),\lambda)d\tau\right]$$

$$\times \sin\frac{s\pi}{p_2 2}$$

$$\left.-\frac{1}{sp_2}\int_{\pi/2}^{\pi} q(\tau)\sin\frac{s}{p_2}(\pi-\tau)w_2(\tau-\Delta(\tau),\lambda)d\tau\right) \qquad (4.14)$$

elde edilir. Burada olası iki durum vardır: 1. $a_2 \neq 0$  2. $a_2 = 0$.

Önce $a_2 \neq 0$ durumunu göz önüne alalım. $\lambda$ yeterince büyük olsun. Eğer (4.14) ü $s$ ile bölersek şunu elde ederiz:

$$-\frac{\gamma_1 a_2}{p_2 \delta_1}\cos\frac{s\pi}{p_1 2}\sin\frac{s\pi}{p_2 2} - \frac{\gamma_2 a_2}{\delta_2 p_1}\sin\frac{s\pi}{p_1 2}\cos\frac{s\pi}{p_2 2} + \frac{ds\gamma_1 a_2}{\delta_1}\cos\frac{s\pi}{p_1 2}\cos\frac{s\pi}{p_2 2}$$

$$-\frac{d\gamma_1 a_1 p_1}{\delta_1}\sin\frac{s\pi}{p_1 2}\cos\frac{s\pi}{p_2 2} - \frac{d\gamma_1}{\delta_1 p_1}\left(\int_0^{\pi/2} q(\tau)\sin\frac{s}{p_1}(\frac{\pi}{2}-\tau)w_1(\tau-\Delta(\tau),\lambda)d\tau\right)\cos\frac{s\pi}{p_2 2}$$

$$-\frac{ds\gamma_2 a_2 p_2}{\delta_2 p_1}\sin\frac{s\pi}{p_1 2}\sin\frac{s\pi}{p_2 2} - \frac{d\gamma_2 p_2 a_1}{\delta_2}\cos\frac{s\pi}{p_1 2}\sin\frac{s\pi}{p_2 2}$$

$$-\frac{d\gamma_2 p_2}{\delta_2 p_1^2}\left(\int_0^{\pi/2} q(\tau)\cos\frac{s}{p_1}(\frac{\pi}{2}-\tau)w_1(\tau-\Delta(\tau),\lambda)d\tau\right)\sin\frac{s\pi}{p_2 2}$$

$$-\frac{d}{p_2}\int_{\pi/2}^{\pi} q(\tau)\sin\frac{s}{p_2}(\pi-\tau)w_2(\tau-\Delta(\tau),\lambda)d\tau = 0.$$

O halde (4.9) ve (4.10) dan (4.14) denklemi yeterince büyük $\lambda$ değerleri için

$$da_2\left(s\frac{\gamma_1}{\delta_1}\cos\frac{s\pi}{p_1 2}\cos\frac{s\pi}{p_2 2} - s\frac{\gamma_2 p_2}{\delta_2 p_1}\sin\frac{s\pi}{p_1 2}\sin\frac{s\pi}{p_2 2}\right) + O(1) = 0$$

şeklinde yeniden yazılabilir.

$\gamma_1\delta_2 p_1 = \gamma_2\delta_1 p_2$ eşitliğini kullanırsak

$$da_2\left(s\frac{\gamma_1}{\delta_1}\cos\frac{s\pi}{p_1 2}\cos\frac{s\pi}{p_2 2} - s\frac{\gamma_1}{\delta_1}\sin\frac{s\pi}{p_1 2}\sin\frac{s\pi}{p_2 2}\right) + O(1) = 0$$

ve buradan da

$$\frac{sda_2\gamma_1}{\delta_1}\cos s\pi\frac{p_1+p_2}{2p_1 p_2} + O(1) = 0. \tag{4.15}$$

Aşikar olarak (4.15) sonsuz sayıda çözüme sahiptir.

$a_2 = 0$ durumunda ise

$$s\sin s\pi\frac{p_1+p_2}{2p_1 p_2} + O(1) = 0$$

halini alır ki denklem sonsuz sayıda çözüme sahiptir. Böylece teorem ispatlanmış oldu. ∎

Şimdi özdeğer ve özfonksiyonların asimptotik özellikleri üzerinde çalışacağız. Bundan sonra $s$ nin yeterince büyük olduğunu varsayacağız.

(4.3) ve (4.9) dan $\left[0,\frac{\pi}{2}\right]$ aralığında

$$w_1(x,\lambda) = O(1). \tag{4.16}$$

(4.4) ve (4.10) dan $\left[\frac{\pi}{2},\pi\right]$ aralığında

$$w_2(x,\lambda) = O(1) \tag{4.17}$$

elde edilir.

$|\lambda| < \infty$ ve $0 \le x \le \frac{\pi}{2}$ için $w_{1s}'(x,\lambda)$; $|\lambda| < \infty$ ve $\frac{\pi}{2} \le x \le \pi$ için $w_{2s}'(x,\lambda)$ türevlerinin varlığı ve sürekliliği Norkin (1972) tarafından gösterilmiştir:

$$w_{1s}'(x,\lambda) = O(1), \quad x \in \left[0,\frac{\pi}{2}\right] \tag{4.18}$$

ve

$$w_{2s}'(x,\lambda) = O(1), \quad x \in \left[\frac{\pi}{2},\pi\right] \tag{4.19}$$

**Teorem 4.3** $n$ bir doğal sayı olsun. Yeterince büyük her $n$ için $a_2 \ne 0$ durumunda (1.1)-

(1.5) probleminin $\dfrac{p_1^2 p_2^2}{(p_1+p_2)^2}(2n+1)^2$ civarında tam olarak bir özdeğeri vardır.

**İspat.** (4.15) denkleminde $O(1)$ ile gösterilen şu ifadeyi göz önüne alacağız:

$$-\dfrac{\gamma_1}{\delta_1}\left\{-da_1 p_1 \sin s\pi\dfrac{p_1+p_2}{2p_1 p_2} + \dfrac{a_2}{p_2}\sin s\pi\dfrac{p_1+p_2}{2p_1 p_2} + \dfrac{a_1\delta_1\gamma_2}{s\delta_2\gamma_1}\cos s\pi\dfrac{p_1+p_2}{2p_1 p_2}\right.$$

$$+\int_0^{\pi/2}\left[\dfrac{\cos s\left(\dfrac{\pi(p_1+p_2)}{2p_1 p_2}-\dfrac{\tau}{p_1}\right)}{sp_1 p_2} + \dfrac{d\sin s\left(\dfrac{\pi(p_1+p_2)}{2p_1 p_2}+\dfrac{\tau}{p_1}\right)}{p_1}\right]q(\tau)w_1(\tau-\Delta(\tau),\lambda)d\tau$$

$$+\left.\int_{\pi/2}^{\pi}\left[\dfrac{\delta_1\cos\dfrac{s}{p_2}(\pi-\tau)}{sp_2^2\gamma_1} + \dfrac{d\delta_1\sin\dfrac{s}{p_2}(\pi-\tau)}{p_2\gamma_1}\right]q(\tau)w_2(\tau-\Delta(\tau),\lambda)d\tau\right\}.$$

(4.16)-(4.19) formülleri göz önüne alındığında gösterilebilir ki büyük $s$ değerleri için bu ifade sonlu türeve sahiptir. Aşikardır ki (4.15) denkleminin kökleri büyük $s$ değerleri için tam sayıların yakınında yer alır. Biz göstereceğiz ki büyük $n$ için (4.15) denkleminin $\dfrac{p_1 p_2}{p_1+p_2}(2n+1)$ civarında sadece bir kökü vardır.

$G(s)=\dfrac{da_2\gamma_1}{\delta_1}s\cos s\pi\dfrac{p_1+p_2}{2p_1 p_2}+O(1)$ fonksiyonunu göz önüne alalım. O halde

$G'(s)=\dfrac{da_2\gamma_1}{\delta_1}\cos s\pi\dfrac{p_1+p_2}{2p_1 p_2}-\dfrac{sda_2\gamma_1\pi(p_1+p_2)}{2p_1 p_2\delta_1}\sin s\pi\dfrac{p_1+p_2}{2p_1 p_2}+O(1)$ türevi $\dfrac{p_1 p_2}{p_1+p_2}(2n+1)$ civarında yeterince büyük $n$ için mevcut olduğundan Rolle teoreminden iddiamızın doğru olduğu görülür. ∎

$n$ yeterince büyük olsun. (1.1)-(1.5) probleminin özdeğerini $\dfrac{p_1^2 p_2^2(2n+1)^2}{(p_1+p_2)^2}$ civarında $\lambda_n=s_n^2$ olarak gösterelim. $s_n=\dfrac{p_1 p_2}{p_1+p_2}(2n+1)+\delta_n$ olarak tanımlayalım. O halde (4.15) den $\delta_n=O(\dfrac{1}{n})$ elde edilir. Sonuç olarak

$$s_n=\dfrac{p_1 p_2}{p_1+p_2}(2n+1)+O(\dfrac{1}{n}). \tag{4.20}$$

(4.20) formülü (1.1)-(1.5) probleminin özfonksiyonunun asimptotik ifadesini elde etmeyi mümkün kılar. (4.3), (4.11) ve (4.16) dan

$$w_1(x,\lambda) = a_2 \cos\frac{s}{p_1}x + O(\frac{1}{s}) \tag{4.21}$$

ve

$$w_1'(x,\lambda) = -\frac{s}{p_1}a_2 \sin\frac{s}{p_1}x + O(1) \tag{4.22}$$

(4.4), (4.17), (4.21), (4.22) ve $\gamma_1\delta_2 p_1 = \gamma_2\delta_1 p_2$ eşitliğinden

$$\begin{aligned}
w_2(x,\lambda) &= \frac{\gamma_1 a_2}{\delta_1}\cos\frac{s\pi}{p_1 2}\cos\frac{s}{p_2}(x-\frac{\pi}{2}) - \frac{\gamma_2 p_2 a_2}{\delta_2 p_1}\sin\frac{s\pi}{p_1 2}\sin\frac{s}{p_2}(x-\frac{\pi}{2}) + O(\frac{1}{s}) \\
&= \frac{\gamma_1 a_2}{\delta_1}\cos\frac{s\pi}{p_1 2}\cos\frac{s}{p_2}(x-\frac{\pi}{2}) - \frac{\gamma_1 a_2}{\delta_1}\sin\frac{s\pi}{p_1 2}\sin\frac{s}{p_2}(x-\frac{\pi}{2}) + O(\frac{1}{s}) \\
&= \frac{\gamma_1 a_2}{\delta_1}\cos(\frac{s\pi}{p_1 2} + \frac{s}{p_2}(x-\frac{\pi}{2})) + O(\frac{1}{s}) \\
&= \frac{\gamma_1 a_2}{\delta_1}\cos s(\frac{\pi(p_2-p_1)}{2p_1 p_2} + \frac{x}{p_2}) + O(\frac{1}{s})
\end{aligned} \tag{4.23}$$

(4.20) yi (4.21) ve (4.23) de yerine koyarsak şunu elde ederiz:

$$u_{1n}(x) = w_1(x,\lambda_n) = a_2 \cos\frac{p_2(2n+1)}{p_1+p_2}x + O(\frac{1}{n}),$$

$$u_{2n}(x) = w_2(x,\lambda_n) = \frac{\gamma_1}{\delta_1}\cos(\frac{\pi(p_2-p_1)(2n+1)}{2(p_1+p_2)} + \frac{p_1(2n+1)}{p_1+p_2}x) + O(\frac{1}{n})$$

Bundan dolayı $u_n(x)$ özfonksiyonları aşağıdaki asimptotik ifadeye sahiptir:

$$u_n(x) = \begin{cases} a_2\cos\dfrac{p_2(2n+1)}{p_1+p_2}x + O(\dfrac{1}{n}), & x\in\left[0,\dfrac{\pi}{2}\right), \\ \dfrac{\gamma_1}{\delta_1}\cos(\dfrac{\pi(p_2-p_1)(2n+1)}{2(p_1+p_2)} + \dfrac{p_1(2n+1)}{p_1+p_2}x) + O(\dfrac{1}{n}), & x\in\left(\dfrac{\pi}{2},\pi\right]. \end{cases}$$

Bazı ek koşullar altında geç kalan argümente bağlı daha kesin asimptotik formüller elde edilebilir. Varsayalım ki aşağıdaki koşullar sağlanır:

a.) $q'(x)$ ve $\Delta''(x)$ türevleri mevcut, $\left[0,\dfrac{\pi}{2}\right)\bigcup\left(\dfrac{\pi}{2},\pi\right]$ aralığında sınırlı ve sırasıyla

$$q'(\frac{\pi}{2}\pm 0) = \lim_{x\to\frac{\pi}{2}\pm 0}q'(x) \text{ ve } \Delta''(\frac{\pi}{2}\pm 0) = \lim_{x\to\frac{\pi}{2}\pm 0}\Delta''(x) \text{ sonlu limitleri mevcuttur.}$$

b.) $\left[0,\dfrac{\pi}{2}\right)\bigcup\left(\dfrac{\pi}{2},\pi\right]$ aralığında $\Delta'(x) \leq 1$, $\Delta(0) = 0$ ve $\lim\limits_{x\to\frac{\pi}{2}+0}\Delta(x) = 0$

b.) yi kullanarak;

$$x - \Delta(x) \geq 0, \qquad x \in \left[0, \frac{\pi}{2}\right) \tag{4.24}$$

ve

$$x - \Delta(x) \geq \frac{\pi}{2}, \qquad x \in \left(\frac{\pi}{2}, \pi\right] \tag{4.25}$$

(4.21), (4.23),(4.24) ve (4.25) den $\left[0, \frac{\pi}{2}\right)$ ve $\left(\frac{\pi}{2}, \pi\right]$ için sırasıyla

$$w_1(\tau - \Delta(\tau), \lambda) = a_2 \cos \frac{s}{p_1}(\tau - \Delta(\tau)) + O(\frac{1}{s}), \tag{4.26}$$

$$w_2(\tau - \Delta(\tau), \lambda) = \frac{\gamma_1 a_2}{\delta_1} \cos s(\frac{\pi(p_2 - p_1)}{2 p_1 p_2} + \frac{\tau - \Delta(\tau)}{p_2}) + O(\frac{1}{s}) \tag{4.27}$$

bulunur. Bu ifadeleri (4.14) de yerine koyarsak;

$$\frac{da_2 \gamma_1 s}{\delta_1} \cos s\pi \frac{p_1 + p_2}{2 p_1 p_2} - \frac{\gamma_1 da_1 p_1}{\delta_1} \sin s\pi \frac{p_1 + p_2}{2 p_1 p_2} - \frac{\gamma_1 a_2}{\delta_1 p_2} \sin s\pi \frac{p_1 + p_2}{2 p_1 p_2} - \frac{a_1 \gamma_2}{s \delta_2} \cos s\pi \frac{p_1 + p_2}{2 p_1 p_2}$$

$$- \frac{\gamma_1}{s \delta_1 p_1 p_2} \int_0^{\pi/2} q(\tau) \cos s(\pi \frac{p_1 + p_2}{2 p_1 p_2} - \frac{\tau}{p_1}) \left[ a_2 \cos \frac{s}{p_1}(\tau - \Delta(\tau)) + O(\frac{1}{s}) \right] d\tau$$

$$- \frac{d \gamma_1}{\delta_1 p_1} \int_0^{\pi/2} q(\tau) \sin s(\pi \frac{p_1 + p_2}{2 p_1 p_2} - \frac{\tau}{p_1}) \left[ a_2 \cos \frac{s}{p_1}(\tau - \Delta(\tau)) + O(\frac{1}{s}) \right] d\tau$$

$$- \frac{1}{s p_2^2} \int_{\pi/2}^{\pi} q(\tau) \cos \frac{s}{p_2}(\pi - \tau) \left[ \frac{\gamma_1 a_2}{\delta_1} \cos s(\pi \frac{p_2 - p_1}{2 p_1 p_2} + \frac{\tau - \Delta(\tau)}{p_2}) + O(\frac{1}{s}) \right] d\tau$$

$$- \frac{d}{p_2} \int_{\pi/2}^{\pi} q(\tau) \sin \frac{s}{p_2}(\pi - \tau) \left[ \frac{\gamma_1 a_2}{\delta_1} \cos s(\pi \frac{p_2 - p_1}{2 p_1 p_2} + \frac{\tau - \Delta(\tau)}{p_2}) + O(\frac{1}{s}) \right] d\tau$$

$$= \frac{da_2 \gamma_1 s}{\delta_1} \cos s\pi \frac{p_1 + p_2}{2 p_1 p_2} - \frac{\gamma_1 da_1 p_1}{\delta_1} \sin s\pi \frac{p_1 + p_2}{2 p_1 p_2} - \frac{\gamma_1 a_2}{\delta_1 p_2} \sin s\pi \frac{p_1 + p_2}{2 p_1 p_2}$$

$$- \frac{d \gamma_1 a_2}{\delta_1} \left[ \frac{1}{p_1} \sin s\pi \frac{p_1 + p_2}{2 p_1 p_2} \int_0^{\pi/2} q(\tau) \cos \frac{s\tau}{p_1} \cos \frac{s}{p_1}(\tau - \Delta(\tau)) d\tau \right.$$

$$- \frac{1}{p_1} \cos s\pi \frac{p_1 + p_2}{2 p_1 p_2} \int_0^{\pi/2} q(\tau) \sin \frac{s\tau}{p_1} \cos \frac{s}{p_1}(\tau - \Delta(\tau)) d\tau$$

$$+ \frac{1}{p_2} \cos s\pi \frac{p_2 - p_1}{2 p_1 p_2} \int_{\pi/2}^{\pi} q(\tau) \sin \frac{s}{p_2}(\pi - \tau) \cos \frac{s}{p_2}(\tau - \Delta(\tau)) d\tau$$

$$+\frac{1}{p_2}\sin s\pi \frac{p_2-p_1}{2p_1p_2}\int_{\pi/2}^{\pi} q(\tau)\sin\frac{s}{p_2}(\pi-\tau)\sin\frac{s}{p_2}(\tau-\Delta(\tau))d\tau\Bigg]+O(\frac{1}{s})$$

$$=\frac{da_2\gamma_1 s}{\delta_1}\cos s\pi \frac{p_1+p_2}{2p_1p_2}-\frac{\gamma_1 da_1 p_1}{\delta_1}\sin s\pi \frac{p_1+p_2}{2p_1p_2}-\frac{\gamma_1 a_2}{\delta_1 p_2}\sin s\pi \frac{p_1+p_2}{2p_1p_2}$$

$$-\frac{d\gamma_1 a_2}{\delta_1}\Bigg[\frac{1}{p_1}\sin s\pi \frac{p_1+p_2}{2p_1p_2}\int_0^{\pi/2} q(\tau)\cos\frac{s\tau}{p_1}\cos\frac{s}{p_1}(\tau-\Delta(\tau))d\tau$$

$$-\frac{1}{p_1}\cos s\pi \frac{p_1+p_2}{2p_1p_2}\int_0^{\pi/2} q(\tau)\sin\frac{s\tau}{p_1}\cos\frac{s}{p_1}(\tau-\Delta(\tau))d\tau$$

$$+\frac{1}{p_2}\cos s\pi \frac{p_2-p_1}{2p_1p_2}\sin\frac{s\pi}{p_2}\int_{\pi/2}^{\pi} q(\tau)\cos\frac{s\tau}{p_2}\cos\frac{s}{p_2}(\tau-\Delta(\tau))d\tau$$

$$-\frac{1}{p_2}\cos s\pi \frac{p_2-p_1}{2p_1p_2}\cos\frac{s\pi}{p_2}\int_{\pi/2}^{\pi} q(\tau)\sin\frac{s\tau}{p_2}\cos\frac{s}{p_2}(\tau-\Delta(\tau))d\tau$$

$$+\frac{1}{p_2}\sin s\pi \frac{p_2-p_1}{2p_1p_2}\sin\frac{s\pi}{p_2}\int_{\pi/2}^{\pi} q(\tau)\cos\frac{s\tau}{p_2}\sin\frac{s}{p_2}(\tau-\Delta(\tau))d\tau$$

$$-\frac{1}{p_2}\sin s\pi \frac{p_2-p_1}{2p_1p_2}\cos\frac{s\pi}{p_2}\int_{\pi/2}^{\pi} q(\tau)\sin\frac{s\tau}{p_2}\sin\frac{s}{p_2}(\tau-\Delta(\tau))d\tau\Bigg]+O(\frac{1}{s})$$

$$=\frac{sda_2\gamma_1}{\delta_1}\cos s\pi \frac{p_1+p_2}{2p_1p_2}-\frac{\gamma_1}{\delta_1}\left(da_1 p_1+\frac{a_2}{p_2}\right)\sin s\pi \frac{p_1+p_2}{2p_1p_2}$$

$$-\frac{d\gamma_1 a_2}{\delta_1}\Bigg[\frac{1}{p_1}\sin s\pi \frac{p_1+p_2}{2p_1p_2}\int_0^{\pi/2} q(\tau)\frac{1}{2}\Bigg[\cos\frac{s\Delta(\tau)}{p_1}+\cos\frac{s}{p_1}(2\tau-\Delta(\tau))\Bigg]d\tau$$

$$-\frac{1}{p_1}\cos s\pi \frac{p_1+p_2}{2p_1p_2}\int_0^{\pi/2} q(\tau)\frac{1}{2}\Bigg[\sin\frac{s\Delta(\tau)}{p_1}+\sin\frac{s}{p_1}(2\tau-\Delta(\tau))\Bigg]d\tau$$

$$+\frac{1}{p_2}\cos s\pi \frac{p_2-p_1}{2p_1p_2}\sin\frac{s\pi}{p_2}\int_{\pi/2}^{\pi} q(\tau)\frac{1}{2}\Bigg[\cos\frac{s\Delta(\tau)}{p_2}+\cos\frac{s}{p_2}(2\tau-\Delta(\tau))\Bigg]d\tau$$

$$-\frac{1}{p_2}\cos s\pi \frac{p_2-p_1}{2p_1p_2}\cos\frac{s\pi}{p_2}\int_{\pi/2}^{\pi} q(\tau)\frac{1}{2}\Bigg[\sin\frac{s\Delta(\tau)}{p_2}+\sin\frac{s}{p_2}(2\tau-\Delta(\tau))\Bigg]d\tau$$

$$-\frac{1}{p_2}\sin s\pi \frac{p_2-p_1}{2p_1p_2}\sin\frac{s\pi}{2}\int_{\pi/2}^{\pi} q(\tau)\frac{1}{2}\Bigg[\sin\frac{s\Delta(\tau)}{p_2}-\sin\frac{s}{p_2}(2\tau-\Delta(\tau))\Bigg]d\tau$$

$$-\frac{1}{p_2}\sin s\pi \frac{p_2-p_1}{2p_1p_2}\cos\frac{s\pi}{p_2}\int_{\pi/2}^{\pi} q(\tau)\frac{1}{2}\Bigg[\cos\frac{s\Delta(\tau)}{p_2}-\cos\frac{s}{p_2}(2\tau-\Delta(\tau))\Bigg]d\tau\Bigg]+O(\frac{1}{s})=0.$$

(4.28)

Şimdi bazı fonksiyonlar tanımlayalım:

$$A(x,s,\Delta(\tau)) = \frac{1}{2}\int_0^x q(\tau)\sin\frac{s}{p_1}\Delta(\tau)d\tau,$$

$$B(x,s,\Delta(\tau)) = \frac{1}{2}\int_0^x q(\tau)\cos\frac{s}{p_1}\Delta(\tau)d\tau.$$

(4.29)

Açıktır ki bu fonksiyonlar $0 < s < \infty$ için $0 \leq x \leq \frac{\pi}{2}$ aralığında sınırlıdır.

$$C(x,s,\Delta(\tau)) = \frac{1}{2}\int_{\pi/2}^x q(\tau)\sin\frac{s}{p_2}\Delta(\tau)d\tau,$$

$$D(x,s,\Delta(\tau)) = \frac{1}{2}\int_{\pi/2}^x q(\tau)\cos\frac{s}{p_2}\Delta(\tau)d\tau.$$

(4.30)

Açıktır ki bu fonksiyonlar $0 < s < \infty$ için $\frac{\pi}{2} \leq x \leq \pi$ aralığında sınırlıdır.

a.) ve b.) koşulları altında aşağıdaki formüller doğrudur (Norkin 1972).

$$\int_0^x q(\tau)\cos\frac{s}{p_1}(2\tau - \Delta(\tau))d\tau = O(\frac{1}{s}),$$

$$\int_0^x q(\tau)\sin\frac{s}{p_1}(2\tau - \Delta(\tau))d\tau = O(\frac{1}{s}),$$

$$\int_{\pi/2}^x q(\tau)\cos\frac{s}{p_2}(2\tau - \Delta(\tau))d\tau = O(\frac{1}{s}),$$

$$\int_{\pi/2}^x q(\tau)\sin\frac{s}{p_2}(2\tau - \Delta(\tau))d\tau = O(\frac{1}{s}).$$

(4.31)

(4.28), (4.29), (4.30) ve (4.31) den şunu elde ederiz:

$$\frac{sd\gamma_1 a_2}{\delta_1}\cos s\pi\frac{p_1+p_2}{2p_1p_2} - \frac{da_1\gamma_1 p_1}{\delta_1}\sin s\pi\frac{p_1+p_2}{2p_1p_2} - \frac{\gamma_1 a_2}{\delta_1 p_2}\sin s\pi\frac{p_1+p_2}{2p_1p_2}$$

$$-\frac{d\gamma_1 a_2}{\delta_1 p_1}B(\frac{\pi}{2},s,\Delta(\tau))\sin s\pi\frac{p_1+p_2}{2p_1p_2} + \frac{d\gamma_1 a_2}{\delta_1 p_1}A(\frac{\pi}{2},s,\Delta(\tau))\cos s\pi\frac{p_1+p_2}{2p_1p_2}$$

$$-\frac{d\gamma_1 a_2}{\delta_1 p_2}D(\pi,s,\Delta(\tau))\sin\frac{s\pi}{p_2}\cos s\pi\frac{p_2-p_1}{2p_1p_2} - \frac{d\gamma_1 a_2}{\delta_1 p_2}C(\pi,s,\Delta(\tau))\sin\frac{s\pi}{p_2}\sin s\pi\frac{p_2-p_1}{2p_1p_2}$$

$$+\frac{d\gamma_1 a_2}{\delta_1 p_2}C(\pi,s,\Delta(\tau))\cos\frac{s\pi}{p_2}\cos s\pi\frac{p_2-p_1}{2p_1p_2}$$

$$-\frac{d\gamma_1 a_2}{\delta_1 p_2}D(\pi,s,\Delta(\tau))\cos\frac{s\pi}{p_2}\sin s\pi\frac{p_2-p_1}{2p_1p_2} + O(\frac{1}{s}) = 0.$$

Eğer bu denklemi $s$ ile bölersek denklem şu hale gelir:

$$\frac{d\gamma_1 a_2}{\delta_1}\cos s\pi\frac{p_1+p_2}{2p_1p_2} - \frac{da_1\gamma_1 p_1}{s\delta_1}\sin s\pi\frac{p_1+p_2}{2p_1p_2} - \frac{\gamma_1 a_2}{s\delta_1 p_2}\sin s\pi\frac{p_1+p_2}{2p_1p_2}$$

$$-\frac{d\gamma_1 a_2}{s\delta_1 p_1}B(\frac{\pi}{2},s,\Delta(\tau))\sin s\pi\frac{p_1+p_2}{2p_1p_2} + \frac{d\gamma_1 a_2}{s\delta_1 p_1}A(\frac{\pi}{2},s,\Delta(\tau))\cos s\pi\frac{p_1+p_2}{2p_1p_2}$$

$$+\frac{d\gamma_1 a_2}{s\delta_1 p_2}C(\pi,s,\Delta(\tau))\cos(\frac{s\pi}{p_2}+s\pi\frac{p_2-p_1}{2p_1p_2})$$

$$-\frac{d\gamma_1 a_2}{s\delta_1 p_2}D(\pi,s,\Delta(\tau))\sin(\frac{s\pi}{p_2}+s\pi\frac{p_2-p_1}{2p_1p_2}) + O(\frac{1}{s^2}) = 0. \tag{4.32}$$

(4.32) denklemini $\sin(\frac{s\pi}{p_2}+s\pi\frac{p_2-p_1}{2p_1p_2}) = \sin s\pi\frac{p_1+p_2}{2p_1p_2}$ ile bölersek

$$\frac{d\gamma_1 a_2}{\delta_1}\cot s\pi\frac{p_1+p_2}{2p_1p_2} - \frac{da_1\gamma_1 p_1}{s\delta_1} - \frac{\gamma_1 a_2}{s\delta_1 p_2}$$

$$-\frac{d\gamma_1 a_2}{s\delta_1 p_1}B(\frac{\pi}{2},s,\Delta(\tau)) + \frac{d\gamma_1 a_2}{s\delta_1 p_1}A(\frac{\pi}{2},s,\Delta(\tau))\cot s\pi\frac{p_1+p_2}{2p_1p_2}$$

$$+\frac{d\gamma_1 a_2}{s\delta_1 p_2}C(\pi,s,\Delta(\tau))\cot s\pi\frac{p_1+p_2}{2p_1p_2} - \frac{d\gamma_1 a_2}{s\delta_1 p_2}D(\pi,s,\Delta(\tau)) + O(\frac{1}{s^2}) = 0$$

olur. Böylece denklem

$$\cot s\pi\frac{p_1+p_2}{2p_1p_2}\left[\frac{d}{sp_2}C(\pi,s,\Delta(\tau)) + d + \frac{d}{sp_1}A(\frac{\pi}{2},s,\Delta(\tau))\right]$$

$$= \frac{1}{s}\left[\frac{d}{p_2}D(\pi,s,\Delta(\tau)) + \frac{da_1 p_1}{a_2} + \frac{1}{p_2} + \frac{d}{p_1}B(\frac{\pi}{2},s,\Delta(\tau))\right] + O(\frac{1}{s^2})$$

ve

$$\cot s\pi\frac{p_1+p_2}{2p_1p_2} = \frac{1}{s}\left[\frac{d}{p_2}D(\pi,s,\Delta(\tau)) + \frac{da_1 p_1}{a_2} + \frac{1}{p_2} + \frac{d}{p_1}B(\frac{\pi}{2},s,\Delta(\tau))\right] + O(\frac{1}{s^2})$$

Eğer $s_n = \frac{p_1 p_2 (2n+1)}{p_1+p_2} + \delta_n$ olarak alırsak

$$\cot(\frac{\pi}{2}(2n+1) + \frac{\pi(p_1+p_2)\delta_n}{2p_1p_2}) = -\tan\frac{\pi(p_1+p_2)\delta_n}{2p_1p_2}$$

$$= \frac{p_1+p_2}{(2n+1)p_1p_2}\left[\frac{d}{p_2}D(\pi,\frac{p_1p_2(2n+1)}{p_1+p_2},\Delta(\tau)) + \frac{da_1 p_1}{a_2} + \frac{1}{p_2} + \frac{d}{p_1}B(\frac{\pi}{2},\frac{p_1p_2(2n+1)}{p_1+p_2},\Delta(\tau))\right]$$

$$+ O(\frac{1}{n^2})$$

olur. Bu yüzden büyük $n$ için

$$\delta_n = -\frac{2}{\pi(2n+1)}\left[\frac{d}{dp_2}D(\pi,\frac{p_1p_2(2n+1)}{p_1+p_2},\Delta(\tau))+\frac{da_1p_1}{a_2}+\frac{1}{p_2}+\frac{d}{dp_1}B(\frac{\pi}{2},\frac{p_1p_2(2n+1)}{p_1+p_2},\Delta(\tau))\right]$$
$$+O(\frac{1}{n^2})$$

olarak bulunur. Sonuç olarak,

$$s_n = \frac{p_1p_2(2n+1)}{p_1+p_2} - \frac{2}{\pi(2n+1)}\left[\frac{d}{dp_2}D(\pi,\frac{p_1p_2(2n+1)}{p_1+p_2},\Delta(\tau))+\frac{da_1p_1}{a_2}+\frac{1}{p_2}\right.$$
$$\left.+\frac{d}{dp_1}B(\frac{\pi}{2},\frac{p_1p_2(2n+1)}{p_1+p_2},\Delta(\tau))\right]+O(\frac{1}{n^2}). \quad (4.33)$$

Böylece sıradaki teoremi ispatlamış olduk.

**Teorem 4.3** Eğer a.) ve b.) koşulları sağlanırsa (1.1)-(1.5) probleminin $\lambda_n = s_n^2$ pozitif özdeğerleri $n \to \infty$ için (4.33) asimptotik formülüne sahiptir.

Özfonksiyonlar için daha kesin bir asimptotik formül elde edebiliriz. (4.3) ve (4.26) dan:

$$w_1(x,\lambda) = a_2\cos\frac{s}{p_1}x - \frac{a_1p_1}{s}\sin\frac{s}{p_1}x - \frac{a_2}{sp_1}\int_0^x q(\tau)\sin\frac{s}{p_1}(x-\tau)\cos\frac{s}{p_1}(\tau-\Delta(\tau))d\tau + O(\frac{1}{s^2}).$$

Böylece (4.29), (4.30) ve (4.31) den

$$w_1(x,\lambda) = a_2\cos\frac{s}{p_1}x - \frac{a_1p_1}{s}\sin\frac{s}{p_1}x$$
$$-\frac{a_2}{sp_1}\int_0^x q(\tau)\left[\sin\frac{s}{p_1}x\cos\frac{s}{p_1}\tau - \cos\frac{s}{p_1}x\sin\frac{s}{p_1}\tau\right]\cos\frac{s}{p_1}(\tau-\Delta(\tau))d\tau + O(\frac{1}{s^2}).$$

$$w_1(x,\lambda) = a_2\cos\frac{s}{p_1}x - \frac{a_1p_1}{s}\sin\frac{s}{p_1}x - \frac{a_2}{sp_1}\int_0^x q(\tau)\left[\sin\frac{s}{p_1}x\cos\frac{s}{p_1}\tau\cos\frac{s}{p_1}(\tau-\Delta(\tau))\right.$$
$$\left.-\cos\frac{s}{p_1}x\sin\frac{s}{p_1}\tau\cos\frac{s}{p_1}(\tau-\Delta(\tau))\right]d\tau + O(\frac{1}{s^2}).$$

$$w_1(x,\lambda) = a_2\cos\frac{s}{p_1}x - \frac{a_1p_1}{s}\sin\frac{s}{p_1}x - \frac{a_2\sin\frac{s}{p_1}x}{2sp_1}\int_0^x q(\tau)\left[\cos\frac{s}{p_1}\Delta(\tau) + \cos\frac{s}{p_1}(2\tau-\Delta(\tau))\right]d\tau$$
$$-\frac{a_2\cos\frac{s}{p_1}x}{2sp_1}\int_0^x q(\tau)\left[\sin\frac{s}{p_1}\Delta(\tau) + \sin\frac{s}{p_1}(2\tau-\Delta(\tau))\right]d\tau + O(\frac{1}{s^2}).$$

$$w_1(x,\lambda) = a_2 \cos\frac{s}{p_1}x - \frac{a_1 p_1}{s}\sin\frac{s}{p_1}x - \frac{a_2 \sin\frac{s}{p_1}x}{sp_1}B(\frac{\pi}{2},s,\Delta(\tau)) + \frac{a_2 \cos\frac{s}{p_1}x}{sp_1}A(\frac{\pi}{2},s,\Delta(\tau))$$
$$+ O(\frac{1}{s^2}).$$

$$w_1(x,\lambda) = a_2 \cos\frac{s}{p_1}x\left[1+\frac{A(\frac{\pi}{2},s,\Delta(\tau))}{sp_1}\right]x - \frac{\sin\frac{s}{p_1}x}{s}\left[a_1 p_1 + \frac{a_2}{p_1}B(\frac{\pi}{2},s,\Delta(\tau))\right] + O(\frac{1}{s^2}). \quad (4.34)$$

$s$ yerine $s_n$ koyarak ve (4.33) ü kullanarak

$$u_{1n} = w_1(x,\lambda_n) = a_2 \cos\frac{p_2(2n+1)}{p_1+p_2}x\left[1+\frac{(p_1+p_2)}{p_1 p_2(2n+1)}A(x,\frac{p_1 p_2(2n+1)}{p_1+p_2},\Delta(\tau))\right]$$
$$+ a_2 \sin\frac{p_2(2n+1)}{p_1+p_2}x\left[\frac{2x}{\pi(2n+1)p_1}\left[\frac{d}{p_2}D(\pi,\frac{p_1 p_2(2n+1)}{p_1+p_2},\Delta(\tau))\right.\right.$$
$$\left.\left.+\frac{da_1 p_1}{a_2}+\frac{1}{p_2}+dB(\frac{\pi}{2},\frac{p_1 p_2(2n+1)}{p_1+p_2},\Delta(\tau))\right]\right] - \frac{p_1+p_2}{p_1 p_2(2n+1)}\sin\frac{p_2(2n+1)}{p_1+p_2}x \quad (4.35)$$
$$\times\left[a_1 p_1 + \frac{a_2}{p_1}B(x,\frac{p_1 p_2(2n+1)}{p_1+p_2},\Delta(\tau))\right] + O(\frac{1}{n^2})$$

elde edilir.

$$\frac{w_1'(x,\lambda)}{s} = -\frac{a_2}{p_1}\sin\frac{s}{p_1}x - \frac{a_1}{s}\cos\frac{s}{p_1}x$$
$$-\frac{a_2}{sp_1^2}\int_0^x q(\tau)\cos\frac{s}{p_1}(x-\tau)\cos\frac{s}{p_1}(\tau-\Delta(\tau)) + O(\frac{1}{s^2})$$
$$= -\frac{a_2}{p_1}\sin\frac{s}{p_1}x - \frac{a_1}{s}\cos\frac{s}{p_1}x$$
$$-\frac{a_2}{sp_1^2}\int_0^x q(\tau)\frac{1}{2}\left[\cos\frac{s}{p_1}(x-\Delta(\tau))+\cos\frac{s}{p_1}(x-(2\tau-\Delta(\tau)))\right] \quad (4.36)$$
$$= -\frac{a_2}{p_1}\sin\frac{s}{p_1}x\left(1+\frac{1}{sp_1}A(x,s,\Delta(\tau))\right)$$
$$-\frac{\cos\frac{s}{p_1}x}{s}\left(a_1 + \frac{a_2}{p_1^2}B(x,s,\Delta(\tau))\right) + O(\frac{1}{s^2}), \qquad x\in\left(0,\frac{\pi}{2}\right].$$

(4.4), (4.27), (4.31), (4.34), (4.36) ve $\gamma_1\delta_2 p_1 = \gamma_2\delta_1 p_2$ eşitliğinden

$$w_2(x,\lambda) = \frac{\gamma_1}{\delta_1}\left\{a_2 \cos\frac{s\pi}{p_1 2}\left[1+\frac{1}{sp_1}A(\frac{\pi}{2},s,\Delta(\tau))\right]-\frac{\sin\frac{s\pi}{p_1 2}}{s}\right.$$

$$\left.\times\left[a_1 p_1+\frac{a_2}{p_1}B(\frac{\pi}{2},s,\Delta(\tau))\right]+O(\frac{1}{s^2})\right\}\cos\frac{s}{p_2}(x-\frac{\pi}{2})$$

$$-\frac{\gamma_2 p_2}{\delta_2 p_1}\left\{a_2 \sin\frac{s\pi}{p_1 2}\left[1+\frac{1}{sp_1}A(\frac{\pi}{2},s,\Delta(\tau))\right]+\frac{\cos\frac{s\pi}{p_1 2}}{s}\right.$$

$$\left.\times\left[a_1 p_1+\frac{a_2}{p_1}B(\frac{\pi}{2},s,\Delta(\tau))\right]+O(\frac{1}{s^2})\right\}\sin\frac{s}{p_2}(x-\frac{\pi}{2})$$

$$-\frac{1}{sp_2}\int_{\pi/2}^{x}q(\tau)\sin\frac{s}{p_2}(x-\tau)\left[\frac{\gamma_1 a_2}{\delta_1}\cos\frac{s}{p_2}(\frac{\pi(p_2-p_1)}{2p_1}+\tau-\Delta(\tau))+O(\frac{1}{s})\right]d\tau$$

$$= \frac{\gamma_1 a_2}{\delta_1}\cos\frac{s\pi}{p_1 2}\cos\frac{s\pi}{p_2 2}\cos\frac{s}{p_2}x\left[1+\frac{1}{sp_1}A(\frac{\pi}{2},s,\Delta(\tau))\right]$$

$$-\frac{\frac{\gamma_1 a_2}{\delta_1}\sin\frac{s\pi}{p_1 2}\cos\frac{s\pi}{p_2 2}\cos\frac{s}{p_2}x}{s}\left[\frac{a_1 p_1}{a_2}+\frac{1}{p_1}B(\frac{\pi}{2},s,\Delta(\tau))\right]$$

$$+\frac{\gamma_1 a_2}{\delta_1}\cos\frac{s\pi}{p_1 2}\sin\frac{s\pi}{p_2 2}\sin\frac{s}{p_2}x\left[1+\frac{1}{sp_1}A(\frac{\pi}{2},s,\Delta(\tau))\right]$$

$$+\frac{\frac{\gamma_1 a_2}{\delta_1}\sin\frac{s\pi}{p_1 2}\sin\frac{s\pi}{p_2 2}\sin\frac{s}{p_2}x}{s}\left[\frac{a_1 p_1}{a_2}+\frac{1}{p_1}B(\frac{\pi}{2},s,\Delta(\tau))\right]$$

$$+\frac{\gamma_1 a_2}{\delta_1}\sin\frac{s\pi}{p_1 2}\sin\frac{s\pi}{p_2 2}\cos\frac{s}{p_2}x\left[1+\frac{1}{sp_1}A(\frac{\pi}{2},s,\Delta(\tau))\right]$$

$$-\frac{\frac{\gamma_1 a_2}{\delta_1}\cos\frac{s\pi}{p_1 2}\sin\frac{s\pi}{p_2 2}\cos\frac{s}{p_2}x}{s}\left[\frac{a_1 p_1}{a_2}+\frac{1}{p_1}B(\frac{\pi}{2},s,\Delta(\tau))\right]$$

$$-\frac{\gamma_1 a_2}{\delta_1}\sin\frac{s\pi}{p_1 2}\cos\frac{s\pi}{p_2 2}\sin\frac{s}{p_2}x\left[1+\frac{1}{sp_1}A(\frac{\pi}{2},s,\Delta(\tau))\right]$$

$$-\frac{\frac{\gamma_1 a_2}{\delta_1}\cos\frac{s\pi}{p_1 2}\cos\frac{s\pi}{p_2 2}\sin\frac{s}{p_2}x}{s}\left[\frac{a_1 p_1}{a_2}+\frac{1}{p_1}B(\frac{\pi}{2},s,\Delta(\tau))\right]$$

$$-\frac{\gamma_1 a_2}{sp_2 \delta_1}\int_{\pi/2}^{x}q(\tau)\frac{1}{2}\left[\sin\frac{s}{p_2}(x+\frac{\pi(p_2-p_1)}{2p_1}-\Delta(\tau))\right.$$

$$\left.+\sin\frac{s}{p_2}(x-\frac{\pi(p_2-p_1)}{2p_1}-(2\tau-\Delta(\tau)))\right]d\tau+O(\frac{1}{s^2})$$

$$= \frac{\gamma_1 a_2}{\delta_1} \cos\frac{s\pi}{p_1 2}\cos\frac{s\pi}{p_2 2}\cos\frac{s}{p_2}x\left[1+\frac{1}{sp_1}A(\frac{\pi}{2},s,\Delta(\tau))\right]$$

$$-\frac{\frac{\gamma_1 a_2}{\delta_1}\sin\frac{s\pi}{p_1 2}\cos\frac{s\pi}{p_2 2}\cos\frac{s}{p_2}x}{s}\left[\frac{a_1 p_1}{a_2}+\frac{1}{p_1}B(\frac{\pi}{2},s,\Delta(\tau))\right]$$

$$+\frac{\gamma_1 a_2}{\delta_1}\cos\frac{s\pi}{p_1 2}\sin\frac{s\pi}{p_2 2}\sin\frac{s}{p_2}x\left[1+\frac{1}{sp_1}A(\frac{\pi}{2},s,\Delta(\tau))\right]$$

$$+\frac{\frac{\gamma_1 a_2}{\delta_1}\sin\frac{s\pi}{p_1 2}\sin\frac{s\pi}{p_2 2}\sin\frac{s}{p_2}x}{s}\left[\frac{a_1 p_1}{a_2}+\frac{1}{p_1}B(\frac{\pi}{2},s,\Delta(\tau))\right]$$

$$+\frac{\gamma_1 a_2}{\delta_1}\sin\frac{s\pi}{p_1 2}\sin\frac{s\pi}{p_2 2}\cos\frac{s}{p_2}x\left[1+\frac{1}{sp_1}A(\frac{\pi}{2},s,\Delta(\tau))\right]$$

$$-\frac{\frac{\gamma_1 a_2}{\delta_1}\cos\frac{s\pi}{p_1 2}\sin\frac{s\pi}{p_2 2}\cos\frac{s}{p_2}x}{s}\left[\frac{a_1 p_1}{a_2}+\frac{1}{p_1}B(\frac{\pi}{2},s,\Delta(\tau))\right]$$

$$-\frac{\gamma_1 a_2}{\delta_1}\sin\frac{s\pi}{p_1 2}\cos\frac{s\pi}{p_2 2}\sin\frac{s}{p_2}x\left[1+\frac{1}{sp_1}A(\frac{\pi}{2},s,\Delta(\tau))\right]$$

$$-\frac{\frac{\gamma_1 a_2}{\delta_1}\cos\frac{s\pi}{p_1 2}\cos\frac{s\pi}{p_2 2}\sin\frac{s}{p_2}x}{s}\left[\frac{a_1 p_1}{a_2}+\frac{1}{p_1}B(\frac{\pi}{2},s,\Delta(\tau))\right]$$

$$-\frac{\gamma_1 a_2}{sp_2\delta_1}\sin\frac{s}{p_2}(x+\frac{\pi(p_2-p_1)}{2p_1})\underbrace{\frac{1}{2}\int_{\pi/2}^{x}q(\tau)\cos\frac{s\Delta(\tau)}{p_2}d\tau}_{D(x,s,\Delta(\tau))}$$

$$+\frac{\gamma_1 a_2}{sp_2\delta_1}\cos\frac{s}{p_2}(x+\frac{\pi(p_2-p_1)}{2p_1})\underbrace{\frac{1}{2}\int_{\pi/2}^{x}q(\tau)\sin\frac{s\Delta(\tau)}{p_2}d\tau}_{C(x,s,\Delta(\tau))}$$

$$\underbrace{-\frac{\gamma_1 a_2}{sp_2\delta_1}\sin\frac{s}{p_2}(x-\frac{\pi(p_2-p_1)}{2p_1})\underbrace{\frac{1}{2}\int_{\pi/2}^{x}q(\tau)\cos(2\tau-\Delta(\tau))d\tau}_{O(\frac{1}{s})}}_{O(\frac{1}{s^2})}$$

$$\underbrace{+\frac{\gamma_1 a_2}{sp_2\delta_1}\cos\frac{s}{p_2}(x-\frac{\pi(p_2-p_1)}{2p_1})\underbrace{\frac{1}{2}\int_{\pi/2}^{x}q(\tau)\sin(2\tau-\Delta(\tau))d\tau}_{O(\frac{1}{s})}}_{O(\frac{1}{s^2})}$$

$$=\frac{\gamma_1 a_2}{\delta_1}\left[1+\frac{1}{sp_1}A(\frac{\pi}{2},s,\Delta(\tau))\right]\left[\cos\frac{s\pi}{p_1 2}\cos\frac{s\pi}{p_2 2}\cos\frac{s}{p_2}x+\cos\frac{s\pi}{p_1 2}\sin\frac{s\pi}{p_2 2}\sin\frac{s}{p_2}x\right.$$

$$+\sin\frac{s\pi}{p_1 2}\sin\frac{s\pi}{p_2 2}\cos\frac{s}{p_2}x - \sin\frac{s\pi}{p_1 2}\cos\frac{s\pi}{p_2 2}\sin\frac{s}{p_2}x\Big] - \frac{\gamma_1 a_2}{s\delta_1 p_2}\Big[\Big[\frac{a_1 p_1 p_2}{a_2} + \frac{p_2}{p_1}B(\frac{\pi}{2},s,\Delta(\tau))\Big]$$

$$\times \Big[\sin\frac{s\pi}{p_1 2}\cos\frac{s\pi}{p_2 2}\cos\frac{s}{p_2}x + \sin\frac{s\pi}{p_1 2}\sin\frac{s\pi}{p_2 2}\sin\frac{s}{p_2}x - \cos\frac{s\pi}{p_1 2}\sin\frac{s\pi}{p_2 2}\cos\frac{s}{p_2}x$$

$$+\cos\frac{s\pi}{p_1 2}\cos\frac{s\pi}{p_2 2}\sin\frac{s}{p_2}x\Big] + \Big[\sin\frac{s}{p_2}(x+\frac{\pi(p_2-p_1)}{2p_1})D(x,s,\Delta(\tau))$$

$$-\cos\frac{s}{p_2}(x+\frac{\pi(p_2-p_1)}{2p_1})C(x,s,\Delta(\tau))\Big]\Big] + O(\frac{1}{s^2})$$

$$= \Big\{\frac{\gamma_1 a_2}{\delta_1}\Big[1+\frac{1}{sp_1}A(\frac{\pi}{2},s,\Delta(\tau))\Big] + \frac{\gamma_1 a_2}{sp_2\delta_1}C(x,s,\Delta(\tau))\Big\}\cos\frac{s}{p_2}(x+\frac{\pi(p_2-p_1)}{2p_1})$$

$$-\Big\{\frac{\gamma_1 a_2}{sp_2\delta_1}D(x,s,\Delta(\tau))+\frac{\gamma_1}{s\delta_1}\Big[a_1 p_1 + \frac{a_2}{p_1}B(\frac{\pi}{2},s,\Delta(\tau))\Big]\Big\}\sin\frac{s}{p_2}(x+\frac{\pi(p_2-p_1)}{2p_1})$$

$$+O(\frac{1}{s^2})\quad , x\in\left(\frac{\pi}{2},\pi\right].$$

$s$ yerine $s_n$ koyarak ve (4.33) eşitliğini kullanarak $x\in\left(\frac{\pi}{2},\pi\right]$ için

$$\begin{aligned}u_{2n} = w_2(x,\lambda_n) &= \Big\{\frac{\gamma_1 a_2}{\delta_1}\Big[1+\frac{p_1+p_2}{p_1^2 p_2(2n+1)}A(\frac{\pi}{2},\frac{p_1 p_2(2n+1)}{p_1+p_2},\Delta(\tau))\Big]\\
&+\frac{\gamma_1(p_1+p_2)}{\delta_1 p_1 p_2^2(2n+1)}C(x,\frac{p_1 p_2(2n+1)}{p_1+p_2},\Delta(\tau))\Big\}\cos(\frac{p_1 x}{p_1+p_2}+\frac{\pi(p_2-p_1)(2n+1)}{2(p_1+p_2)})\\
&+\Big\{\frac{\gamma_1 a_2}{\delta_1}\Big[\frac{2}{\pi(2n+1)}\Big[\frac{d}{p_2}D(\pi,\frac{p_1 p_2(2n+1)}{p_1+p_2},\Delta(\tau))+\frac{da_1 p_1}{a_2}+\frac{1}{p_2}\\
&+\frac{d}{p_1}B(\frac{\pi}{2},\frac{p_1 p_2(2n+1)}{p_1+p_2},\Delta(\tau))\Big]\Big]\Big(\frac{x}{p_2}+\frac{\pi(p_2-p_1)}{2p_1 p_2}\Big)\\
&-\Big[\frac{a_2\gamma_1(p_1+p_2)}{p_1 p_2^2\delta_1(2n+1)}D(x,\frac{p_1 p_2(2n+1)}{p_1+p_2},\Delta(\tau))+\frac{\gamma_1(p_1+p_2)}{\delta_1 p_1 p_2(2n+1)}\big(a_1 p_1\\
&+\frac{a_2}{p_1}B(\frac{\pi}{2},\frac{p_1 p_2(2n+1)}{p_1+p_2},\Delta(\tau))\big)\Big]\Big\}\sin(\frac{p_1 x}{p_1+p_2}+\frac{\pi(p_2-p_1)(2n+1)}{2(p_1+p_2)})+O(\frac{1}{n^2}).\end{aligned}\quad(4.37)$$

$a_2 = 0$ durumu için de benzer işlemler yapılarak özdeğer ve özfonksiyonlar için asimptotik formüller elde edilebilir.

## 5. SONUÇ

Son olarak elde edilen bulgular sonucunda aşağıdaki teorem ispatlanmış oldu:

**Teorem 5.1**   Eğer a.) ve b.) koşulları sağlanıyorsa (1.1)-(1.5) probleminin $u_n(x)$ özfonksiyonlarının $n \to \infty$ için asimptotik gösterilimi

$$u_n(x) = \begin{cases} u_{1n}(x), & x \in [0, \frac{\pi}{2}), \\ u_{2n}(x), & x \in (\frac{\pi}{2}, \pi] \end{cases}$$

şeklindedir. Burada $u_{1n}(x)$, (4.35) ve $u_{2n}(x)$, (4.37) formülleriyle tanımlıdır.

## 6. KAYNAKLAR

# ÖZGEÇMİŞ

Erdoğan Şen, 28 Aralık 1985 tarihinde Kırcaali'de doğmuştur. 2004 yılında Pertevniyal Anadolu Lisesini bitirdikten sonra aynı yıl Gebze Yüksek Teknoloji Enstitüsü Fen Fakültesi Matematik Bölümü'nde lisans eğitimine başlamıştır. Lisans Eğitimini 2008 yılında tamamladıktan sonra 2009 yılında Yıldız Teknik Üniversitesi Fen Bilimleri Enstitüsü Matematik Anabilim Dalı'nda yüksek lisans eğitimine başlamıştır. Tez aşamasında ise Namık Kemal Üniversitesi Fen Bilimleri Enstitüsü Matematik Anabilim Dalı'na yatay geçiş yapmıştır. 2009 yılının Aralık ayından beri Namık Kemal Üniversitesi Fen Edebiyat Fakültesi Matematik Bölümü'nde araştırma görevlisi olarak çalışmaktadır.